\documentclass[12pt]{article}
\usepackage{amsmath,amsfonts,amssymb,amsthm,color}

\newcommand{\I}{{\bf 1}}
\newtheorem{proposition}{Proposition}[section]
\newtheorem{theorem}[proposition]{Theorem}

\newtheorem{remark}[proposition]{Remark}

\newtheorem{example}[proposition]{Example}

\newcommand{\nc}{\newcommand}
\nc{\R}{{\mathbb R}}
\nc{\N}{{\mathbb N}}
\nc{\Z}{{\mathbb Z}}

\nc{\BP}{\mathbb{P}}
\nc{\BE}{\mathbb{E}}
\nc{\BY}{\mathbb{Y}}
\nc{\BQ}{\mathbb{Q}}
\nc{\cB}{\mathcal{B}}
\nc{\bN}{{\mathbf N}}

\begin{document}

\author{K.\ Borovkov\footnote{
Department of Mathematics and Statistics, University of Melbourne,
k.borovkov@ms.unimelb.edu.au.}
\ and \ G. Last
\footnote{
Institut f\"ur Stochastik, Karlsruher Institut f\"ur  Technologie,
guenter.last@kit.edu.} }
\title{On Rice's formula for stationary multivariate\\
piecewise smooth  processes}
\date{\today}

\maketitle

\begin{abstract}
\noindent
Let $X=\{X_t: t\ge 0\}$  be a stationary piecewise continuous $\R^d$-valued
process that moves between jumps along the integral curves of a
given continuous vector
field, $S\subset\R^d$ be a smooth surface. The aim of this paper is to derive a
multivariate version of Rice's formula, relating the intensity of the
point process of (localized) continuous crossings  of $S$ by $X$ to
the distribution of $X_0$.
Our result is illustrated by examples relating to queueing networks
and stress release network models.
\end{abstract}

\noindent
{\em Keywords:}  level crossings, Rice's formula, stationarity,
Palm probabilities, piecewise-deterministic process, stochastic network

\vspace{0.1cm}
\noindent
2000 Mathematics Subject Classification: Primary 60J75; secondary 60G55.

\section{Introduction}
\setcounter{equation}{0}

The classical Rice's formula going back to~\cite{Rice}
gives the intensity $\nu (u)$ of
crossings (originally, upcrossings) of a given level $u$
by a univariate continuous
stationary Gaussian process $X_t$ in terms of the joint distribution
of $(X_t,X'_t)\stackrel{d}{=}(X_0, X'_0)$,
the process' value and its derivative at a fixed time
(provided that the derivative exists in some suitable sense,
e.g.\ in mean quadratic):
\begin{align} \label{R-1}
\nu (u) =\int  |z|\, p(u,z)\, dz,
\end{align}
where $p (\cdot, \cdot)$ is the joint density
$(X_0, X'_0)$ which is assumed to exist.
Later on the result has been extended to more general
classes of differentiable (in some
suitable sense) stationary processes, covering not only the first
moments but also
higher order factorial moments of the numbers of crossing,
and even to more general
settings for continuous random processes and fields.
The formula proved to be quite
useful in a number of applied areas, including signal processing,
reliability, sea waves
and others. For detailed accounts of the history of results of this
kind and further
bibliography, the interested reader is referred to \cite{Ry00}, \cite{LeSp04} and
Chapter~3 in~\cite{AzWs09}.

The case of processes with jumps and degenerate finite-dimensional
distributions drew
much less attention, although, from the applications' viewpoint,
it is scarcely less
interesting than the one of continuous  processes. However,
the heuristics behind the
formula based on ``Kac's counting formula" giving the
number of crossings of a level $u$
by a $C^1$-function $f$ on $[0,1]$ as
\begin{align}\label{Kac}
\lim_{\delta\to 0+} \frac{1}{2\delta}
\int_0^1 |f' (t)|\, \I \{|f(t) - u| <\delta\} \,dt
\end{align}
(under a couple of further technical assumptions and denoting
by $\I A$ the indicator
function of the set $A$),   seems to be applicable in that
case as well, provided that
the process jumps at finite intensity and is smooth between the jump times.
We note that \eqref{Kac} is a consequence
of Federer's coarea theorem  (see e.g.\ (7.4.15) in~\cite{AT07})
\begin{align}\label{coarea}
\int_0^1 g(t)\,|f' (t)| \,dt
=\int^\infty_{-\infty}\sum_{s\in[0,1]} g(s)
\,\I\{f(s)=t\}\,dt,
\end{align}
applied to the function $g(s):=\I\{|f(s) - u| <\delta\}$.

An analogue of \eqref{R-1} for the intensity $\nu_c (u)$ of
continuous level crossings by
general univariate piecewise deterministic Markov processes that has the form
$$
\nu_c (u) = |\mu (u)| p(u),
$$
where $\mu (\cdot)$ is the drift coefficient
of the process and $p(\cdot)$  the density
of $X_0$, was established in~\cite{BoLa08} (see also Theorem~\ref{triceone}
below; one should mention here an earlier paper~\cite{BaNe72}
where the case of Poisson shot-noise
processes was considered). The result was used in~\cite{BoLa08}
to obtain the asymptotic
behaviour of the point processes of high level crossings
(i.e.\ as $u\to\infty$) in a
number of interesting and important for applications special cases.

The proof in  \cite{BoLa08} relied on the Markov structure of
the process and in fact did not assume the existence of
the density $p$\,---\,its existence was part of the
assertion of the main theorem there. The natural question on
whether Rice's formula for
piecewise smooth processes can be extended to the
multivariate and non-Markovian cases
remained open. In the present paper we give a positive answer to it.

The paper is organized as follows. In Section~2, we  describe the main class of
processes we will be working with and then present the main result
together with its
proof. Section~3 presents examples to illustrate our main result.

\section{The main result}
\label{sec2}
\setcounter{equation}{0}

First we will describe the  main model of multivariate random processes
$X=\{X_t: t\ge
0\}$ dealt with in this paper. The two key elements of the model
are a point process
$N=\{N(B): B\in \cB (\R_+)\}$ of jumps in our process $X$
(here and in what follows,
$\cB (\cdot)$ denotes  the class of Borel subsets of $\cdot$)
and a vector field $\mu:D
\to  \R^d$ defined on an open domain $D\subset \R^d$ and
specifying the dynamics of~$X$
between the jumps. Note that we allow trivial jumps (of zero size) as well.

We assume that the following assumptions are satisfied.

\medskip\hangindent=0.7cm \hangafter=0\noindent
{\bf (A.1)}~$N$ is a stationary simple counting
process on $\R_+$, which has a finite
intensity $\lambda_N:=\BE N((0,1])$ and is such that $N(\R_+)=\infty$ a.s.

\medskip

The latter  implies that the process $N$ is locally finite
and hence its points can be
enumerated in the increasing order. We denote them by
$0 < T_1 <T_2 <\dots$   and set
$T_0:=0$ for convenience (this is not a point of $N$ a.s.).

\medskip\hangindent=0.7cm \hangafter=0\noindent
{\bf (A.2)}~$\mu\in C^1 (D)$.
\medskip

This assumption implies that there exist continuous functions
$t_-:\R^d \to  (-\infty,0)$ and $t_+:\R^d \to  (0,\infty)$
such that, for any $x\in D$, there exists a unique $C^1$-function
$q(x,\cdot):(t_-(x),t_+(x)) \to  D$
satisfying the integral equation
\begin{align}\label{ie}
q(x,t)=x+\int^t_0 \mu(q(x,s))\, ds,\quad t\in (t_-(x),t_+(x))
\end{align}
(Picard-Lindel\"of theorem, see e.g.\ p.8  in~\cite{Hartman}).
Moreover, for any  fixed
$x\in \R^d$ there is a neighbourhood of $(x,0)\in\R^{d+1}$
in which $q(\cdot, \cdot)$
will also be continuously differentiable (Peano's theorem on
dependence on initial conditions, see e.g.\ p.95 in~\cite{Hartman}).

The integral curves $q$ specify the dynamics of the process $X$ between its jumps.

\medskip\hangindent=0.7cm \hangafter=0\noindent
{\bf (A.3)}~Assume that, for any $n\ge 0,$ one has
$X_{T_n}\in D$, $T_{n+1} <  T_n + t_+(X_{T_n})$   and
\begin{align}\label{piece}
X_t=q(X_{T_n},t-T_n), \quad T_n\le t<T_{n+1}.
\end{align}
Moreover, $X$ and $N$ are  jointly stationary, i.e.\ the
distribution of the bivariate process $\{ (X_{s+t},N((s,s+t])): t\ge 0\}$ does
not depend on $s\ge 0$.

\medskip

Next we will list assumptions involving the surface $S$ of which the continuous
crossings by $X$ we are concerned with. The latter are defined as follows:
we say that
$X$ has a {\em continuous crossing\/} of $S$ at time $s>0$
if $X_{s-}=X_s\in S$ and
there is a $\delta>0$ such that $X_t\notin S$ for $t\in
(s-\delta,s+\delta)\setminus\{s\}$.

\medskip\hangindent=0.7cm \hangafter=0\noindent
{\bf (A.4)}~Let $S\subset D$ be the relative interior of a
$(d-1)$-dimensional (not
necessarily connected) $C^1$-manifold with or without boundary, and
$\{n(x):x\in S\}$
be a continuous field of unit normals  to $S$. Denoting by
$\langle \cdot ,\cdot \rangle$ the Euclidean scalar product in $\R^d$,
we assume that
\begin{align} \label{pos}
\langle n(x),\mu(x) \rangle\ne 0,\quad x\in S.
\end{align}

\begin{remark}\label{r1a}\rm
Let $\tau_x:=\inf\{t>0:q(x,t)\in S\}$ be the first positive
time the integral curve of
$\mu$ leaving  from $x$ at time zero hits the surface~$S$.
It is not hard to see that,
if $S'$ is a compact subset of $S$, then from {\bf (A.2)} and {\bf (A.4)}
it follows that
\begin{align} \label{2.3}
\inf\{\tau_x: x\in S'\}>0.
\end{align}
This (together with the fact that $X$ jumps only finitely often in finite time
intervals) implies that the times of continuous crossings of
$S$ through a compact subset of $S$ cannot accumulate
in finite time.
\end{remark}

The times of continuous crossings of $S$ by $X$ form an
at most countable set $N_c$ that will be identified with
a random counting measure on $[0,\infty)$.
Then
\begin{align}\label{2.3a}
\Phi_c(C):=\sum_{s\in N_c}\I\{(s,X_s)\in C\}
=\int \I\{(s,X_s)\in C\}N_c(ds),
 \quad C\in\cB \bigl([0,\infty)\times\R^d\bigr),
\end{align}
defines a random (integer-valued) measure $\Phi_c$ on
$[0,\infty)\times\R^d$. For $t\ge 0$ and
$S'\in \cB(\R^d)$, the random variable $\Phi_c([0,t]\times S')$
need not be finite. However, if $S'$ is a compact subset
of $S$ then \eqref{2.3} implies that
$\Phi_c([0,t]\times S')<\infty$. Moreover, since $N$
has a finite intensity,
\begin{align}\label{2.3b}
\nu_c(B):=\BE \Phi_c((0,1]\times B),\quad B\in\cB (\R^d),
\end{align}
is finite, whenever $B$ is a compact subset of $S$. Therefore $\nu_c(\cdot)$ is a
$\sigma$-finite measure on $\cB(\R^d)$. For any compact 
$B\subset S$, the point process
$\Phi_c(\cdot\times B)$ is stationary. This is enough to 
derive the (refined) Campbell theorem stating that
\begin{align} \label{refC}
\BE \int g(s,X_s)\,N_c(ds)=\iint g(s,x)\,ds\,\nu_c(dx).
\end{align}
for any measurable function $g: \R_+\times\R^d \to  \R_+$,
cf.\ e.g.\ (1.2.19) in~\cite{BB}.

\begin{remark}\label{r2}\rm
Assuming that $A$ is a small enough open set to ensure 
that $\nu_c(S\cap A)<\infty$,
observe that $\nu_c(S\cap A)^{-1}\nu_c(\cdot)$ can be 
interpreted as the distribution of
the value of $X$ at a {\em typical\/} time of  continuous 
crossing of~$S\cap A$. This is
a particular instance of a  Palm distribution, see e.g.~\cite{BB}.
\end{remark}

The Palm measure $\pi_0$ of the pairs of values of $X$ 
just before and after a typical jump of $X$ is defined by
\begin{align}\label{palm-}
\pi_0(B):=\BE\sum_{n=1}^\infty
\I\{T_n\le 1,X_{T_n-}\ne X_{T_n},(X_{T_n-},X_{T_n})\in B\},
\quad B\in\cB(\R^d\times\R^d).
\end{align}
Note that $\pi_0(\R^d\times\R^d)\le \lambda_N<\infty$.

\medskip\hangindent=0.7cm \hangafter=0\noindent
{\bf (A.5)}~The distribution $\pi$ of $X_0$ has a continuous density $p$ in a
neighborhood of $S$,  and
\begin{align}\label{asspalm}
\min \bigl\{ \pi_0((\R^d\setminus S)\times S) ,\,
  \pi_0(S \times(\R^d\setminus S))\bigr\}=0.
\end{align}

Now we are ready to state our main result.

\begin{theorem} \label{trice}
Under assumptions  {\bf (A.1)}\hspace{.1 mm}--\hspace{.1 mm}{\bf (A.5)},
one has the identity
\begin{align} \label{rice4}
\nu_c(B)=\int_{S\cap B}
|\langle n(x),\mu(x)\rangle|\, p(x)\, \mathcal{H}^{d-1}(dx),
\quad B\in\cB (\R^d),
\end{align}
where $\mathcal{H}^{d-1}$ is the $(d-1)$-dimensional Hausdorff
measure on~$\R^d$.
\end{theorem}

\begin{remark}\label{r87}\rm As it will be seen from the first half
of the proof of Theorem~\ref{trice}, the full continuity 
assumption on  $p$  (which is
part of~{\bf (A.5)}) can actually be somewhat weakened to the 
boundedness of $p$   in a
neighborhood of $S$ and its right-continuity 
(in case the first number in \eqref{asspalm} vanishes) 
on $S$ along the flow meaning that
$p(q(x,0+))=p(x)$, $x\in S$.
\end{remark}

In the one-dimensional case the above theorem simplifies to the 
following assertion.

\begin{theorem} \label{triceone}
In the case  $d=1$, assuming that $S=\{u\}$ for some $u\in D$
such that $\mu (u) \neq 0$, and that\/
{\bf (A.1)}\hspace{.1 mm}--\hspace{.1 mm}{\bf (A.3)} and\/ {\bf (A.5)}
are satisfied, one has
\begin{align}\label{rice1}
\nu_c (\{u\})=|\mu(u)| p(u).
\end{align}
\end{theorem}

\begin{remark}\label{r87a}\rm
In the Markovian case, representation \eqref{rice1} was established 
in \cite{BoLa08}.
More precisely, it was shown there that there exists a density $p$
satisfying~\eqref{rice1}. Due to the Markovian structure of the process, it was
possible to derive the result under weaker technical assumptions.
\end{remark}

In the case when
\[
S=S^u:=\{x\in \R^d:x_1=u\}
\]
for some $u\in\R$, a continuous crossing of $S$ is a continuous 
crossing of the level
$u$ by the first component of $X$. In this case Theorem 
\ref{trice} takes the following form.

\begin{theorem} \label{trice3}
Let assumptions  {\bf (A.1)}\hspace{.1 mm}--\hspace{.1 mm}{\bf (A.3)} 
be satisfied and
$u\in\R$ be such that $S^u\subset D$ and $\mu_1(x)\ne 0$ for all 
$x\in S^u$, where
$\mu_1$ is the first component of $\mu$. Assume that {\bf (A.5)} holds 
with $S=S^u$. Then, for $B\in\cB(\R^{d}),$
\begin{align*}
\nu_{c}(B)=\idotsint \I_B(u,x_2,\ldots,x_d)\, |\mu_1(u,x_2,\ldots,x_d)|\,
p(u,x_2,\ldots,x_d)\, dx_2 \cdots dx_d.
\end{align*}
\end{theorem}

\begin{remark}\label{r93a}\rm
Theorem \ref{trice3} is another, more straigthforward, generalization of 
\eqref{rice1}.
Assume now that $0<\nu_{c}(S^u)<\infty$  and consider a 
``typical time'' of a continuous
crossing of the level $u$ by the first component of $X$. 
Then the measure $\BQ_u (\cdot) :=\nu_{c}(S^u)^{-1}\nu_c(\{u\}\times\cdot)$ 
describes the distribution of the other
components of $X$ at this time. This distribution can be 
interpreted in terms of the
drift-modulated density $p_1$ proportional to $|\mu_1(x)|p(x)$ 
(assuming that $\BE|\mu_1(X_0)|<\infty$). 
If $(Y_1,\ldots,Y_d)$ is a random vector with density $p_1$,
then $\BQ_u$ is the conditional distribution of $(Y_2,\ldots,Y_d)$ 
given that $Y_1=u$.
\end{remark}

\begin{remark}\label{r95}\rm Let $k\in\{1,\ldots,d\}$ and
assume that $S=\tilde{S}\times\R^{d-k}$, where $\tilde{S}\subset\R^k$
is a $(k-1)$-dimensional smooth surface.
Let $\{\tilde{n}(x):x\in \tilde{S}\}$
be a continuous field of unit normals  to $\tilde{S}$. Let
$\tilde{X}:=(X^{(1)},\ldots,X^{(k)})$ 
and $Y:=(X^{(k+1)},\ldots,X^{(d)})$, where $X=(X^{(1)},\ldots,X^{(d)})$.
There is a one-to-one correspondence between the
continuous crossings of $S$ by the process $X$ and
continuous crossings of $\tilde{S}$ by the process $\tilde{X}$.
Equation \eqref{rice4} can be written as
\begin{align}\label{rice6}
\nu_c(B)=\int_\mathbb{Y}\int _{\tilde{S}}
|\langle \tilde{n}(x),\tilde\mu(x,y)\rangle|
\I_B(x,y)p(x|y)\mathcal{H}^{k-1}(dx)
\BP(Y_0\in dy),
\quad B\in\mathcal{B}(\R^k\times\BY),
\end{align}
where $\BY:=\R^{d-k}$, $\tilde\mu(x,y)$ is the vector of the first
$k$ components of $\mu(x,y)$, and $x\mapsto p(x|y)$ is the
conditional density of $\tilde{X}_0$
given that $Y_0=y$. In this form the result might
be generalizable to other stationary
pairs $(\tilde{X},Y)$. The process $\tilde{X}$ should
remain piecewise deterministic for given $Y$.
But the process $Y$ might take values in a more general
space $\BY$. 
In this paper we will make no attempt to establish
such an extension of our results.
\end{remark}

To prove Theorem~\ref{trice}, we will need an auxiliary
result that  requires some
further  notation. First of all, for our purposes it will
suffice  that that result
would hold in a ``local setting", i.e.\ for $S$ replaced with
$S\cap A,$ where $A$ is a
small enough open subset of~$\R^d$. As can easily be seen from
the observation that we
made after stating assumption {\bf (A.2)} and from {\bf (A.4)},
if we understand by $S$ such a ``small piece" of the original surface,
then  the following will be satisfied:

\medskip\hangindent=0.7cm \hangafter=0\noindent
{\bf (A.6)}~The surface $S$ is connected and
relatively compact, \eqref{2.3} holds with
$S'=S$ and $\nu_c:= \nu_c(S) <\infty.$ Furthermore,
there exists a $u_0>0$ such that
$t_+(x)\ge u_0$ for all $x\in S$ and, for any  $u\in [0,u_0]$,
\begin{align*}
S_u:=\{q(x,u):x\in S\}
\end{align*}
is a $C^1$-surface with a continuous field $\{n_u(x):x\in S_u\}$
of unit normals to it satisfying
\begin{align}\label{pos2}
\inf\bigl\{\langle n_u(x),\mu(x) \rangle: \  x\in S_u,\  u\in [0, u_0]\bigr\}> 0.
\end{align}
Moreover, $\pi$ has a density $p$ in a
neighbourhood of\/ $S_{(0,u_0)}$, where
$$
S_{I}:=\bigcup_{u\in I} S_u, \quad I\subset \R.
$$

Now denote by $N^u_{c}$   the stationary point process of the times of
all continuous crossings of $S_u$ by~$X$. For  any
$C\in\cB([0,\infty)\times\R^d),$
let $\Phi^u_c(C)$ be the number of all
$s\in N^u_{c}$ such that $(s,X_s)\in C$ and
$\nu^u_{c}(B):=\BE\Phi^u_c([0,1]\times B)$, $B\in\cB(\R^d)$.

\begin{proposition}\label{p92}
Under assumptions {\bf (A.1)}\hspace{.1 mm}--\hspace{.1 mm}{\bf (A.4)}
and {\bf (A.6)},
for any measurable function  $g:\R^d \to \R_+$, one has
\begin{align}
 \label{67}
\int g(x)\nu_c^u(dx)
 =
\int_{S_u}  |\langle n(x),\mu(x)\rangle |\, g(x)p(x)\mathcal{H}^{d-1}(dx)
\end{align}
for $\mathcal{H}^1$-almost all $u\in[0,u_0]$.
\end{proposition}

\noindent
{\sc Proof.} For any $j\ge 0$ set $T'_j:=T_j\wedge 1$ and,
in particular, $T'_0:=0$. For $j\ge 1$ we define
$$
I_j:=(T'_{j-1},T'_j),
\qquad L_j:=\{X_t:t\in I_j\}=\{q(X_{T'_{j-1}},t-T'_{j-1}):t\in I_j\}.
$$
Fix a $B\in\cB(\R^d)$ and assume that $u\in (0,u_0)$. By definition,
$\Phi^u_c(I_j\times B)>0$ if and only if $L_j\cap S_u\cap B\ne\emptyset$.
On the other
hand, \eqref{pos2} implies   that $\Phi^u_c(I_j\times B)\le 1$, so that
$$
\Phi^u_c(I_j\times B)=\I \{L_j\cap S_u\cap B\ne\emptyset\}.
$$
Therefore
\begin{align*}
\Phi^u_c((0,1)\times B)=\sum^\infty_{j=1}\I\{L_j\cap S_u\cap B\ne\emptyset\}
\end{align*}
and, for any $v\in(0,u_0),$
\begin{align} \label{2.31}
\int^v_0 \Phi^u_c((0,1)\times B)\, du
=\sum^\infty_{j=1}\int^v_0\I\{L_j\cap S_u\cap B\ne\emptyset\}du.
\end{align}
Now set
\[
J_j (v):= \{t\in I_j:  X_t\in S_{(0,v)}\}, \quad
 U_j (v) := \{u\in (0,v): L_j\cap S_u\ne\emptyset\}.
\]
Clearly, the last two sets are either simultaneously empty or are open
intervals of the
same length; in the latter case, put $u_j (v):= \inf U_j (v)$.
Therefore,
\begin{align}
\int^v_0\I\{L_j\cap S_u\cap B\ne\emptyset\}du
 &= \int_{U_j (v)} \I \{q (X_{u_j(v)},  u-u_j (v) )\in B\} du
  \notag\\
 &= \int_{J_j (v)} \I \{X_t\in B\} dt
 = \int_{I_j (v)} \I \{X_t\in S_{(0,v)} \cap B\} dt,
  \label{G*}
\end{align}
so that \eqref{2.31} becomes
$$
\int^v_0 \Phi^u_c((0,1)\times B)\, du
 =\int^1_0\I\{X_t\in S_{(0,v)}\cap B\}dt.
$$
Taking expectations on both sides of the last relation and using
Fubini's  theorem and stationarity of $X$, we obtain that
\begin{align*}
\int^v_0 \nu^u_c(B)du=\BE\int^1_0\I\{X_t\in S_{(0,v)}\cap B\}dt
=\BP(X_0\in S_{(0,v)}\cap B).
\end{align*}
As functions of $B\in \cB (\R^d)$, both sides specify a measure,
and so the standard
argument shows that, for any measurable function $g:\R^d \to \R_+$,
\begin{align} \label{2.31a}
\int_0^v du \int g(x)\nu_c^u(dx)
 =  \int_{S_{(0,v)}}  g(x) p(x) \mathcal{H}^{d}(dx).
\end{align}

Now we can assume without loss of generality that $S$ admits
a $C^1$-parametrization
$(w_1,\ldots,w_{d-1})  \mapsto z(w_1,\ldots,w_{d-1})$,
where $(w_1,\ldots,w_{d-1})$ varies in an open set $W\subset \R^{d-1}$.
For $(w_1,\ldots,w_{d-1})\in W$ and $u\in [0,u_0]$, define
$$
\psi(w_1,\ldots,w_{d-1},u):=q(z(w_1,\ldots,w_{d-1}),u),
$$
which, for a fixed $u\in [0,u_0]$, will be a
$C^1$-parametrization of the ``parallel"  surface $S_u$.

Next we denote by $\partial_i$ the operator of partial differentiation
with respect to
$w_i$, $i=1,\ldots,d-1$, and let $\partial_d \psi:=\partial/\partial u$.
A simple linear
algebra calculations shows that the Jacobian $J\psi$ of
$\psi=\psi(w_1,\ldots,w_{d-1},u)$ satisfies
$$
|J\psi|= |\langle n_u (\psi) , \partial_d\psi\rangle |\, H
 \equiv |\langle n_u (\psi) , \mu (\psi)\rangle |\, H ,
$$
where $H^2 =H^2(w_1,\ldots,w_{d-1},u)$ is the determinant of
the matrix $\bigl(\langle
\partial_i \psi,\partial_j \psi\rangle\bigl)_{i,j=1,\ldots,d-1}$.
However, for any fixed
$u\in[0,u_0]$, $H(w_1,\ldots,w_{d-1},u)\, dw_1\cdots dw_{d-1}$
is the surface element of
$S_u$ in the coordinates $(w_1,\ldots,w_{d-1})$,  so that changing
coordinates on the
right-hand side of \eqref{2.31a} yields
\begin{align}
\int_{S_{(0,v)}}  g(x) p(x) \mathcal{H}^{d}(dx) &=
\int_{W\times (0,v)}   g(\psi) p(\psi)\, |J\psi| \,dw_1 \cdots  dw_{d-1}\, du
 \notag\\
 &= \int_{W\times (0,v)}  g(\psi) p(\psi) \,
  |\langle n_u (\psi) , \mu (\psi)\rangle |\, H \,dw_1 \cdots  dw_{d-1}\, du
  \notag\\
  &=  \int_0^v du  \int_{S_u} g(x) p(x) \,
  |\langle n_u (x) , \mu (x)\rangle |\, \mathcal{H}^{d-1}(dx),
   \label{G**}
\end{align}
which immediately implies the assertion of Proposition~\ref{p92}.\qed

\begin{remark}\label{r87aa}\rm
Assume that $f$ is a real-valued $C^1$-function defined on
an open domain  $\tilde D\subset R^d,$ with non-vanishing gradient
and such that $S_u=\{x\in\tilde D:f(x)=u\}$
for all small enough~$u$. Such a function exists,
at least for suitably small pieces
of~$S$. We may then apply Federer's coarea theorem
(see e.g.\ (7.4.15) in~\cite{AT07})
on each open interval $(T'_{j-1},T'_j)$ to the level sets
of the function $t\mapsto
f(X_t)$. While this would provide an alternative way for
deriving~\eqref{G*}, we have
preferred to give a direct argument presented in the above proposition.
In a quite similar spirit the coarea theorem can be used to
derive Rice's formula for smooth
processes, see Section 11.4 in~\cite{AT07}.
It was actually U.~Z\"ahle who first used
in~\cite{Za84} the coarea theorem to prove Rice's formula for certain continuous
processes. We also note in passing that the coarea formula could be used to
establish~\eqref{G**} as well. However, our more explicit
argument yields additional
information that is needed in the proof of Theorem~\ref{trice}.
\end{remark}

\noindent {\sc Proof of Theorem \ref{trice}.} Since both sides of 
\eqref{rice4} are
$\sigma$-additive in $B$, it is no restriction of generality to 
assume that assumption
{\bf (A.6)} is satisfied.
Moreover, we can assume that $S$ admits a smooth parametrization as in the
proof of Proposition \ref{p92}. This is due to the fact that the surface 
$S$ can be represented as a ``mosaic'' of ``small pieces'' 
for which the assumption will be
satisfied owing to assumptions 
{\bf (A.1)}\hspace{.1 mm}--\hspace{.1mm}{\bf (A.4)} on
the original~$S$.

Furthermore, it is not hard to see that, to prove the theorem, it suffices to
demonstrate that~\eqref{67} holds at $u=0$ for continuous and bounded~$g$. 
We will show
that by proving  
that, under the assumption $\pi_0((\R^d\setminus S)\times
S)=0$, both sides of~\eqref{67} are right-continuous at $u=0$, 
as Proposition~\ref{p92}
will imply then the desired result. The case when 
only the second term on the left-hand
side of \eqref{asspalm} turns into zero 
(i.e.\ $\pi_0(S \times(\R^d\setminus S))=0$) can
be dealt with in exactly the same way by establishing the 
left-continuity of both sides
of~\eqref{67} at $u=0$ in this situation 
(essentially via a time-reversal argument).

Using the notation from the proof of Proposition~\ref{p92} and setting
\[
h_u(w) := \langle n_u(\psi),\mu(\psi)\rangle g(\psi)p(\psi) ,
 \quad w=(w_1,\ldots,w_{d-1}), \quad \psi=\psi(w,u),
\]
we have, for $u\in[0,u_0],$
\begin{align}
 \label{2.67}
\int _{S_u}
\langle n_u(x),\mu(x)\rangle g(x)p(x)\mathcal{H}^{d-1}(dx)
=\int h_u(w) H(w,u) \, \mathcal{H}^{d-1}(dw).
\end{align}
As noted after stating assumption {\bf (A.2)},
one has $q\in C^1$,  and so $n_u(\psi(y,u))$ is a
continuous function of $u$, leading to
$$
h_{0+} (w) =  \langle n(z(w)),\mu(z(w))\rangle g(z(w))p(z(w)).
$$
Similarly, as $u\to 0+$, $H^2(w,u)$ converges to the value of the
determinant of the
matrix
$\bigl( \langle \partial_i z(w),\partial_j z(w)\rangle\bigr)_{i,j=1,\ldots,d-1}$.
Now the dominated  convergence theorem implies   that~\eqref{2.67}
converges to  the right-hand side of \eqref{67} at~$u=0$.

To establish the desired right-continuity of the left-hand side of \eqref{67}
we assume that \eqref{pos2} holds. Introduce the
following point process $\Phi_d$ on $\R_+\times\R^d\times\R^d$:
\begin{align*}
\Phi_d(\cdot)
 :=\sum_{n=1}^\infty \I\{X_{T_n-}\ne X_{T_n}\}
 \I\{(T_n,X_{T_n-},X_{T_n})\in \cdot\}.
\end{align*}
Let $u\in[0,u_0]$ and $t\ge 0$. A continuous crossing  
of $S_u$ can only occur on a
trajectory of $X$ that arrives at the surface from the 
{\em inside\/} of $S_{(0,u)}$
along an integral curve of $\mu$ (cf.~\eqref{pos2}). 
Therefore each such crossing of
$S_u$ should be preceded by an entry to $S_{[0,u)}$, either 
along a drift line or by
jump. Taking into account the possibility of having $X_0\in S_{[0,u)}$, 
we obtain the bound
\begin{align*}
N^u_c([u,t+u]) \le N_c([0,t+u]) + \Phi_d\bigl([0,t+u]\times
(\R^d\setminus S_{[0,u]} )\times S_{[0,u]} \bigr)+1.
\end{align*}
Therefore
\begin{align*}
\int^{t+u}_u g(X_s)N^u_c(ds)
 & \le \int^{t+u}_0(g(X_s)+\varepsilon(u))N_c(ds)\\
 & +{g}^*\bigl[\Phi_d
\bigl([0,t+u]\times (\R^d\setminus S_{[0,u]} )\times S_{[0,u]}\bigr)+1\bigr],
\end{align*}
where ${g}^*:= \sup_x g(x)$ and
\begin{align} \label{2.67a}
\varepsilon(u):=\sup\{|g(q(x,v))-g(x)|:x\in S,0\le v\le u\} \to 0  \quad
\mbox{as}\quad   u\to 0+
\end{align}
due to the uniform continuity of the mapping $(x,u)\mapsto g(q(x,u))$
on $\overline{S}\times[0,u_0]$, $\overline{S}$ denoting the closure of $S$.

Now taking expectations on both sides of the obtained inequality and
using Campbell's formula~\eqref{refC} yields
\begin{align*}
t\int g(x)\nu^u_c(dx) \le
(t+u)\int (g(x)+\varepsilon(u))\nu_c(dx) +
{g}^*\bigl( (t+u)\pi_0((\R^d\setminus S_{[0,u]})\times S_{[0,u]})+1\bigr),
\end{align*}
where we also used Campbell's theorem for 
$\pi_0 (\cdot) =\BE\Phi_d([0,1]\times\cdot)$.
After dividing  by $t$ and letting $t\to\infty$, we obtain
\begin{align*}
\int g(x)\nu^u_c(dx) \le\int g(x)\nu_c(dx)+\varepsilon(u)\nu_c(S) + {g}^*
\pi_0((\R^d\setminus S_{[0,u]})\times S_{[0,u]}).
\end{align*}
In view of \eqref{2.67a} and the fact that the assumption 
$\pi_0((\R^d\setminus S)\times S)=0$ implies that 
$\pi_0((\R^d\setminus S_{[0,u]})\times S_{[0,u]})\to 0$ as $u\to 0+$,
this leads to
\begin{align*}
\limsup_{u\to 0+}\int g(x)\nu^u_c(dx)\le \int g(x)\nu_c(dx).
\end{align*}

To derive the converse inequality, we start with the observation 
that any continuous
crossing of $S$ in $[0,t]$ is followed either by a continuous crossing of
$S_u$  or by a jump from $S_{(0,u]}$ to its complement  within the time interval
$[0,t+u]$, so that
\begin{align*}
N^u_c([0,t+u]) \ge N_c([0,t]) - 
\Phi_d\bigl([0,t+u]\times S_{(0,u]}\times(\R^d\setminus S_{(0,u]})\bigr)-1.
\end{align*}
Next, similarly to our argument above, we obtain
\begin{align*}
\int g(x)\nu^u_c(dx) \ge \int g(x)\nu_c(dx)-\varepsilon(u)\nu_c(S)
   -  {g}^* \pi_0(S_{(0,u]}\times(\R^d\setminus S_{(0,u]})).
\end{align*}
Since $\lim_{u\to 0+} S_{(0,u]} =\varnothing$, it is clear that the
continuity of $\pi_0$ implies now that
\begin{align*}
\liminf_{u\to 0+}\int g(x)\nu^u_c(dx)\ge \int g(x)\nu_c(dx),
\end{align*}
which completes the proof of the theorem.
\qed

\section{Examples}\label{examples}
\setcounter{equation}{0}

In this section we will present two examples showing possible uses
of our main result.

\begin{example} \label{queue}\rm
Consider a general queueing network model with $d$ servers operating in
stationary
regime, with arrivals of customers (possibly in batches) to the
network being governed
by a stationary simple point process. Each customer, upon
completion of its service at
node $j\in\{1,\ldots, d\}$ of the network,   proceeds to another node
for further service or leaves the network, according to some routing
mechanism. All the arrival, transition and departure
times form a stationary point process~$N$, and it is at 
these times that the state of
the process $X_t=(X^{(1)}_t,\ldots,X^{(d)}_t)\in\R^d$ describing the 
residual workloads on
the nodes can change by a jump. Between the events, the values 
of $X_t$ decrease
according to the relation $\frac{d}{dt}X_t = \mu (X_t)$ for some $C^1$-function
$\mu:\R^d \to\R_-^d$, so that the service rate at node $j$ can depend 
on the residual
workload at the node and, moreover, it can even depend on the workloads 
at other nodes
$i\ne j$ as well. To make this description compatible with the assumptions in
Section~2, we allow $X^{(j)}_t<0$ interpreting as the residual 
workload at node $j$ at
time $t$ the value $\max\{X^{(j)}_t, 0\}$, and let $D:=\R^d$.

For $i\in\{1,\ldots,d\}$ let $S_i:=\{x=(x_1,\ldots,x_d)\in\R^d:x_i=0\}$. Then the
continuous crossing of the surface $S_i$ corresponds to server $i$ 
becoming idle. Let
$\nu_i(B)$ denote the intensity of these crossings through a point in $B\in
\cB(\R^{d})$. Provided that the assumptions of Theorem~\ref{trice3} 
are satisfied, we obtain
\begin{align}\label{123}
\nu_{i}(B)=\idotsint \I_B(x^{i})\,|\mu_i(x^{i})|\,p(x^{i})\, dx_1\cdots dx_{i-1}
dx_{i+1}\cdots dx_d,
\end{align}
for $B\in\cB(\R^{d})$, where
$x^i:=(x_1,\ldots, x_{i-1},0,x_{i+1},\ldots,x_d)$ and
$\mu_i$ is the $i$th component of $\mu$.
The normalization of \eqref{123} yields the (Palm)
distribution of the network at a typcial departure
time from node $i$.

Note that assumptions {\bf (A.1)}\hspace{.1 mm}--\hspace{.1 mm}{\bf (A.3)} 
are rather
mild and that Theorem~\ref{trice3} also requires $\mu_i(x)<0$ for $x\in S_i$. 
In assumption \eqref{asspalm} only the condition 
$\pi_0(S_i\times(\R^d\setminus S_i))=0$ is of relevance. 
This assumption says that if there 
is a jump at an instant
when server $i$ becomes empty, then the workload of 
this server is not allowed to
increase by this jump, neither by an internal transition 
(including feedback) nor by an
external arrival. Again, this is a rather weak assumption.


We can also consider the ``composite surface" 
$S:=\bigcup_i S_{(i)}$, where $S_{(i)}$
is the set of all $x\in\R^d$ with $x_i=0$ and 
$x_j\ne 0$ for $j\ne i$. (Under weak assumptions any continuous 
crossing of $S_i$ is also a continuous crossing of
$S_{(i)}$.) Theorem~\ref{trice3} provides the Palm distribution 
of the residual workloads at the time when one of the servers
becomes idle while all the others are
still working. The probability of server $i$ becoming 
idle given a typical instant when
(exactly) one of the servers becomes idle is then given
by  $\nu_c(S_{(i)})/\nu_c(S)$. We skip further details.
\end{example}

\begin{example} \label{stress}\rm
The classical stress release model in seismology (see e.g.\ \cite{BoVe00} and
references to earlier work therein) is a piecewise deterministic 
Markov process $X_t$
representing the level of ``stress'' at a seismic 
fault at time $t$. The value $X_t$
continuously increases at a linear rate due to the 
tectonic loading of the fault and
drops by  random jumps when the stress discharges by way of 
earthquakes that occur at
random times whose intensity is given by $\psi (X_t)$ for some suitably chosen
increasing risk function $\psi$ (e.g.\ $\psi (x) =e^{\beta x}$ 
for some $\beta >0$).
Note that the remote measuring of stress levels at seismic faults is an extremely
difficult problem, so the value $X_t$ is usually not 
observable. All the information on
the process one can have access to is contained in the times, 
locations and magnitudes of jumps.

A more interesting multinode analog of the model was discussed
in~\cite{BoBe03}, where
it was demonstrated, in particular, that already a two-node
stress release network can
reproduce the famous Omori's law for the intensity of earthquake aftershocks.

In the multinode model, the values of the components of the random process
$X_t=(X^{(1)}_t, \ldots,X^{(d)}_t )\in\R^d$ represent the time $t$  stress levels at
individual seismic faults $j\in\{1,\ldots, d\}$ constituting a local fault system.
Between jumps, the dynamics of the process are given by 
$\frac{d}{dt}X_t = \mu$ for a
constant vector $\mu \in\R^d$, Note that one can have $\mu_j <0$ 
which corresponds to
tectonic unloading of   stress at node~$j$ 
(of course, we can consider a more general
model with a variable $\mu$ as well; similar remarks apply 
to all the other elements of
the model construction). Jumps (``seismic events'') 
occurrence at node $j$ is driven by a
Markovian random mechanism with the probability of a jump occurring 
at the node in the
infinitesimal time interval $dt$ given by $\psi_j (X_{t-})dt$ for a given risk
function~$\psi_j(x)$.

When the $n$th seismic event occurs at node $j$ 
(say, at time $T_{j,n}$), the value of
stress at the node changes by a random quantity 
$\xi_{j,n}$, $n=1,2,\ldots,$  which may
be assumed to be i.i.d.\ random variables. Moreover, 
the stress levels at other nodes
can also change at that instance: for a given constant 
transfer matrix $(r_{ij})\in
\R^{d\times d}$, one has 
$ X^{(i)}_{T_{j,n}} = X^{(i)}_{T_{j,n}-} + r_{ij} \xi_{j,n}$, $i\ne j$
(for more detail, see~\cite{BoBe03}).

One of the main problems one hopes to be able to solve 
in mathematical seismology is to
give advanced earthquake warnings. Within the framework of 
the multinode stress release
model, that warning would have to be given at the time 
when the  cumulative jump
intensity $\sum_{j=1}^d \psi (X^{(j)}_t)$ exceeds a given threshold 
$u>0$. That is, we
are looking at continuous crossings of the surface 
$S:=\bigl\{x: \sum_{j=1}^d \psi(x_j) = u\bigr\}$ 
by our process $X_t$. Our main result allows to find the distribution
of $X_t$ at the (typical) time of such crossing and hence, for example, 
to derive the probability
for a given fault to trigger the forthcoming seismic event.
\end{example}

\noindent {\bf Acknowledgements.} This research was supported by a grant 
of the German
Science Foundation (DFG) and the ARC Centre of Excellence for Mathematics and
Statistics of Complex Systems (MASCOS).

\end{document}